\documentclass[a4paper,numbook,babel,final,envcountsame,11pt]{article}

\usepackage{amssymb}
\usepackage{amsthm}
\usepackage{amsmath}
\usepackage{graphicx}
\usepackage{array,hhline}

\newtheorem{lemma}{Lemma}[section]
\newtheorem{theorem}[lemma]{Theorem}
\newtheorem{proposition}[lemma]{Proposition}

\newtheorem{corollary}[lemma]{Corollary}
\theoremstyle{definition}
\newtheorem{definition}[lemma]{Definition}

\numberwithin{equation}{section}
\numberwithin{figure}{section}

\newcommand{\Lset}{\mathcal{L}}
\newcommand{\Mset}{\mathcal{M}}

\newcommand{\Sset}{\mathcal{S}}

\newcommand{\Xset}{\mathcal{X}}

\newcommand{\Zset}{\mathcal{Z}}

\begin{document}

\title{\huge On the Noether and the Cayley-Bacharach theorems with PD multiplicities}         
\author{Hakop Hakopian, Navasard Vardanyan}        
\date{}          

\maketitle

\centerline{\sl Yerevan State University, Yerevan, Armenia}

\centerline{E-mails: \, {\it hakop@ysu.am; \,\, vardanyan.navasard@gmail.com}}

\vskip4mm

\begin{abstract} In this paper we prove the Noether theorem with the multiplicities described by
 PD operators. Despite the known analog versions in this case the provided conditions are 
necessary and sufficient. We also prove the Cayley-Bacharach theorem with PD multiplicities. 
  As far as we know this is the first generalization of this theorem for multiple intersections.
\end{abstract}


\vskip4mm

\noindent{\bf MSC2010 numbers:} 41A05, 41A63, 14H50.

\vskip1mm

\noindent{\bf Keywords:} polynomial interpolation; $n$-independent set; PD multiplicity space; arithmetical multiplicity.

\baselineskip16pt

\section{Introduction}

Let $\Pi$ be the space of all bivariate polynomials. Let also $\Pi_n$ be the space of bivariate polynomials of total degree
at most $n:$
\begin{equation*}
\Pi_n=\left\{\sum_{i+j\leq{n}}a_{ij}x^iy^j
\right\}.
\end{equation*}
We have that
\begin{equation} \label{N=}
N:=\dim \Pi_n=\binom{n+2}{2}.
\end{equation}
Consider a set of $s$ linear operators (functionals) on $\Pi_n:$
\begin{equation*}
{\mathcal L}_s=\{ L_1, \dots ,  L_s \} .
\end{equation*}
The problem of finding a polynomial $p \in \Pi_n$ which satisfies
the conditions
\begin{equation}\label{int cond}
L_ip = c_i, \ \ \quad i = 1, 2, \dots s,
\end{equation}
is called the Lagrange interpolation problem with operators.

In our paper we consider linear operators $L$ which are partial differential operators evaluated at points:
$$L f= p\left(\frac{\partial }{\partial x},\frac{\partial }{\partial y}\right)f\big|_{(x_0,y_0)},$$
where $p\in\Pi.$
We say that $L$ has degree $d,$ where $d=\deg p.$

\begin{definition}
A set of  operators ${\mathcal L}_s$ is called \emph{$n$-correct} if
for any data $\{c_1, \dots, c_s\}$ there exists a unique polynomial
$p \in \Pi_n$, satisfying the conditions \eqref{int cond}.
\end{definition}
A necessary condition of $n$-correctness of $\Lset_s$ is:
$|{\mathcal L}_s|=s= N.$

A polynomial $p \in \Pi_n$ is called an \emph{$n$-fundamental
polynomial} for an operator $L_k\in {\mathcal X}_s$ if
\begin{equation*}
L_ip = \delta _{i k}, \  i = 1, \dots , s ,
\end{equation*}
where $\delta$ is the Kronecker symbol.

We denote the
$n$-fundamental polynomial for $L \in{\mathcal L}_s$ by $p_L^\star=p_{{L}, {\mathcal L}}^\star.$ Sometimes we also call fundamental a polynomial at which vanish all operators but one, since it is a nonzero
constant times the fundamental polynomial.

The following is a Linear Algebra fact:
\begin{proposition} \label{prp:poised}
The set of operators ${\mathcal L}_N$, with $|{\mathcal L}_N|=N=\binom{n+2}{2},$ is $n$-poised if and only if the
following implication holds:
$$p \in \Pi_n\ \hbox{and}\ L_ip = 0, \ i = 1, \dots , N \Rightarrow p = 0.$$
\end{proposition}

\subsection{$n$-independent and $n$-dependent sets}

Next we introduce an important concept of $n$-dependence of sets of operators:
\begin{definition}
A set of operators ${\mathcal L}$ is called \emph{$n$-independent} if
each operator has a fundamental polynomial in $\Pi_n$. Otherwise, ${\mathcal
L}$ is called \emph{$n$-dependent.}
\end{definition}
\noindent Clearly fundamental polynomials are linearly independent.
Therefore a necessary condition of $n$-independence of the set ${\mathcal L}$ is $|{\mathcal
L}| \le N.$

Suppose $\lambda$ is a point in the plane. Consider the operator $L_\lambda$ defined by $L_\lambda f=f(\lambda).$
We say that a set of points $\Xset$ is $n$-independent ($n$-correct) if the set
of operators $\{L_\lambda:\lambda \in\Xset\}$ is $n$-independent ($n$-correct).

Suppose a set of operators ${\mathcal L}$ is $n$-independent. Then by using the Lagrange formula:
\begin{equation*}
p = \sum_{L\in{\mathcal L}} c_{L} p_{{L}, {\mathcal L}}^\star,\quad c_L=Lp,
\end{equation*}
we obtain a polynomial $p \in \Pi_n$ satisfying the
interpolation conditions \eqref{int cond}.

Thus we get a simple characterization of $n$-independence:

\noindent A node set ${\mathcal L}_s$ is $n$-independent if and only
if the interpolation problem \eqref{int cond} is
\emph{$n$-solvable,} meaning that for any data $\{c_1, \dots , c_s
\}$ there exists a (not necessarily unique) polynomial $p \in \Pi_n$
satisfying the conditions \eqref{int cond}.

Now suppose that ${\mathcal L}_s$ is $n$-dependent. Then some operator
$L_{i_0}, \ i_0 \in \{1, \dots , s\}$, does not possess an $n$-fundamental
polynomial. This means that the following implication holds$:$
$$p \in \Pi_n,\ L_{i_0}p = 0\ \forall i \in \{1, \dots , s\} \setminus\{i_0\} \Rightarrow L_{i_0}p =
0.$$

Let $\ell$ be a line.   We say that $p\in\Pi$ vanishes at $\lambda\in \ell$ with the multiplicity
 $m$ if $$(D_a)^ip\big|_\lambda=0,\ i=0,\dots,m-1,$$ where $a||\ell$ and $D_a$ is the directional derivative.

The following proposition is well-known (see, e.g., \cite{HJZ09b}
Proposition 1.3):
\begin{proposition}\label{n+1points}
Suppose that  $\ell$ is a line and a polynomial $p \in
\Pi_n$ vanishes at some points of $\ell$ with the sum of multiplicities $n+1.$ Then we have
\begin{equation} \label{pl}
p = \ell r, \ \text{where} \ r\in\Pi_{n-1}.
\end{equation}
\end{proposition}
Note that this relation also yields that the mentioned $n+1$ conditions are independent, since $\dim\Pi_n-\dim\Pi_{n-1}=n+1.$


\subsection{Multiple intersections}
Let us start with the following well-known
relation for polynomial $R$ and functions $g$ and $f$ (see, e.g., \cite{H03}, formula (16)):
\begin{equation}\label{main}R(D) [gf] =\sum_{i,j\ge 0} {1\over{i!j!}} g^{(i,j)}R^{(i,j)}(D) f. \end{equation}
Here we use the following notations
$$R(D):=R(\frac{\partial }{\partial x},\frac{\partial }{\partial y}),\quad
R^{(i,j)}:=D^{(i,j)}R:=\left(\frac{\partial }{\partial x}\right)^i\left(\frac{\partial }{\partial y}\right)^jR.$$
Notice that to verify \eqref{main} it suffices to check it for R being a monomial, which
reduces \eqref{main} to Leibniz's rule.

To simplify notation, we shall use the same letter $p$,
say, to denote the polynomial $p$ and the curve given by the
equation $p(x,y)=0.$ Thus the notation $\lambda\in p$ means that the point $\lambda$ belongs to the curve $p(x,y)=0.$ Similarly $p\cap q$ for polynomials $p$ and $q$ stands for  the set of intersection points of the curves $p(x,y)=0$ and $q(x,y)=0.$

Below we bring the definition of multiplicities described by PD operators (see \cite{MMM}, \cite{H04}, \cite{HT}):
\begin{definition} The following space is called the multiplicity space of the polynomial $p\in\Pi_n$ at the point $\lambda\in p:$
$$\Mset_\lambda(p)=\left\{h\in \Pi : D^\alpha h(D) p(\lambda)=0\ \forall \alpha\in \mathbb{Z}_+^2\right\}.$$
\end{definition}
\noindent Denote by $\Zset_0=p\cap q$ the set of intersection points of curves (polynomials) $p$ and $q.$
\begin{definition} Suppose that $p,q\in \Pi$ and  $\lambda\in \Zset_0.$ Then the following space is called the multiplicity space of the intersection point $\lambda:$
$$\Mset_\lambda(p,q)=\Mset_\lambda(p)\cap \Mset_\lambda(q).$$
We have that (see \cite{H04}) the spaces $\Mset_\lambda(p,q)$ are $D$-invariant, meaning that
\begin{equation}\label{inv} f\in \Mset_\lambda(p,q) \Rightarrow \frac{\partial f}{\partial x}\ \hbox{and}\ \frac{\partial f}{\partial y}\in \Mset_\lambda(p,q).
\end{equation}
The number $\dim \Mset_\lambda(p,q)$ is called the arithmetical multiplicity of the point $\lambda.$
\end{definition}
Denote
$$\Mset(p,q)=\bigcup_{\lambda\in \Zset_0}\Mset_\lambda(p,q).
$$
 We say that $f\in \Pi_k$ vanishes at ${\mathcal M}_\lambda (p, q)$ if $h(D)f(\lambda) = 0\ \forall h\in {\mathcal M}_\lambda (p, q).$

We say also that the polynomials $p$ and $q$ have no intersection point at infinity if the leading homogeneous parts of $p$ and $q$ have no common factor.

\begin{theorem} [\cite{H04}, Theorem 3]\label{B} Suppose that polynomials $p,q\in \Pi,\ \deg p=m,\ \deg q=n,$ have no intersection point at infinity.
Then the number of the intersection points, counted with the arithmetical multiplicities, equals $mn:$
\begin{equation}\label{abcde}\sum_{\lambda\in \Zset_0}\dim\Mset_\lambda(p,q) =mn.\end{equation}
\end{theorem}
Let us bring the formulation of this result in the homogeneous case.
Let $\Pi_n^0$ be the space of trivariate homogeneous polynomials of total degree
$n.$ In analog way we are defining the multiplicity space $\Mset_\lambda^0(p,q).$
\begin{theorem} [\cite{H04}, Corollary 3] Suppose that polynomials $p\in\Pi_m^0, q\in \Pi_n^0$ have no common component.
Then the number of the intersection points, counted with the arithmetical multiplicities, equals $mn:$
$$\sum_{\lambda\in \Zset_0}\dim\Mset_\lambda^0(p,q) =mn.$$
\end{theorem}
\section {The Noether theorem}

Suppose that $p,q\in \Pi,\ \deg p=m,\ \deg q=n,$ and $p\cap q:=\{\lambda_1,\ldots,\lambda_s\}.$
Let us choose a basis in the space $\Mset_{\lambda_k}(p,q)$ in the following way. Let $\{L_{m1}^{k},\ldots, L_{m{i_{m}}}^{k}\}$ be a maximal independent set of linear  operators with the highest degree $m:=m_k.$
Next we choose $\{L_{m-11}^{k},\ldots, L_{m-1{i_{m-1}}}^{k}\}$ to be a maximal independent set of linear  operators with the degree $m-1.$ Continuing similarly for the degree $0$ we have only one operator $L_{01}^k.$

It is easily seen that the above operators $L_{\mu i}^k,$ form a basis in the linear space $\Mset_{\lambda_k}(p,q).$
Denote
$$\Lset^k(p,q):=\Lset^{\lambda_k}(p,q):=\bigcup_{i,\mu}L_{\mu i}^k,\qquad \Lset(p,q):=\bigcup_{k}\Lset^k(p,q).$$
Notice that, according to Theorem \ref{B}, we have that $|\Lset(p,q)|=mn,$
provided that $p$ and $q$ have no intersection point at infinity.

\begin{lemma} The set of linear operators $\Lset(p,q)$ is $\gamma_0$-independent for sufficiently large $\gamma_0.$
\end{lemma}
\textbf{Proof.}
Consider the set of  the linear operators of fixed node $\lambda_{k_0}=(x_0,y_0)$ of degrees up to $\nu,$ i.e.,
$$\Sset_{\nu,k_0}:=\bigcup_{\mu \le \nu}L_{\mu i}^{k_0}.$$
Let us first  find a fundamental polynomial $p^*$ for an operator of the highest degree $\nu,$ say, for $L_{\nu1}^{k_0}$ within $\Sset_{\nu,k_0}.$
We seek $p^*$ in the form
\begin{equation*}
p^*(x,y)=\sum_{i+j=\nu}a_{ij}(x-x_0)^i(y-y_0)^j.
\end{equation*}
Then we readily get that
$L_{\mu i}^{k_0}p^*=0,\ \hbox{if}\ \mu\le \nu-1.$ Now suppose that
$$
L_{\nu s}^{k_0}f= p_s\left(\frac{\partial }{\partial x},\frac{\partial }{\partial y}\right)f\big|_{(x_0,y_0)},\quad s=1,\ldots,i_s,
$$
where $p_s(x,y)=\sum_{i+j\le\nu}b_{ij}^s(x-x_0)^i(y-y_0)^j.$
Then the conditions of the fundamentality of $p^*$ reduce to the following linear system:
$$
L_{\nu i}^kp^*=\sum_{i+j=\nu}a_{ij}b_{ij}^si!j!=\delta_{ij},\quad s=1,\ldots,i_s.
$$
The linear independence of highest degrees of the operators $L_{\nu i}^k$ means the independence of the vectors  $\{b_{ij}^s\}_{i+j=\nu}.$
Hence the above system has a solution.

Now notice that to complete the proof it is enough to obtain  a fundamental polynomial of $L_{\nu i}^k$ over the set $\Sset_{\nu,k_0}\cup\bigcup_{k\neq k_0}\Lset^k(p,q).$
To this purpose for each $k\in \{1,\ldots,s\}\setminus \{k_0\}$ consider $m_k$ lines passing  through $\lambda_k,$ and not passing through $\lambda_{k_0}.$
Then by multiplying $p^*$ by the product of these lines we obtain, in view of the formula \eqref{main}, a polynomial which is
a desired fundamental polynomial.
\hfill$\Box$

Next, we are going to prove the Noether theorem with the multiplicities described by PD operators.
\begin{theorem} \label{Noether} Suppose that polynomials $p,q\in \Pi,\ \deg p=m,\ \deg q=n,$ have no intersection point at infinity. Suppose also that $f\in \Pi_k$ vanishes at $\Mset_\lambda(p,q)$ for each $\lambda\in p\cap q.$ Then we have that
\begin{equation}\label{aaa}f=Ap+Bq,
\end{equation}
where $A\in \Pi_{k-m},\ B\in \Pi_{k-n}.$
 \end{theorem}
Note that the inverse theorem is true. Indeed, if \eqref{aaa} holds then $f\in \Pi_k$ and, in view of the formula \eqref{main}, we have that  and $f$ vanishes at $\Mset_\lambda(p,q)$ for each $\lambda\in p\cap q.$

\textbf{Proof.}
{\bf Step 1.}
Suppose that $k\ge k_0=\max\{m+n,\gamma_0\},$ where $\gamma_0$ is chosen such that the set of linear operators $\Lset(p,q)$ is $\gamma_0$-independent.

Consider two linear spaces
$${\mathcal V}=\left\{ f\in\Pi_k\ :\ f \ \hbox{vanishes at}\ \Mset_\lambda(p,q)\ \forall \lambda\in p\cap q\right\},$$
$${\mathcal W}=\left\{ Ap+Bq\ :\ A\in \Pi_{k-m},\ B\in \Pi_{k-n}\right\}.$$
In view of the formula \eqref{main} we have that ${\mathcal W}\subset {\mathcal V}.$
To prove the relation \eqref{aaa} we need to verify that ${\mathcal W}= {\mathcal V}.$ To this end it suffices to show that $\dim{\mathcal W}= \dim{\mathcal V}.$

Since the set of linear operators $\Lset(p,q)$ is $\gamma_0$-independent we obtain readily that the set is also $k$-independent, where $k\ge \gamma_0.$

Hence, in view of Theorem \ref{B}, we have that $$\dim{\mathcal V}=\dim\Pi_k- |\Lset(p,q)| = {k+2\choose 2}-mn.$$
Denote
$${\mathcal W_1}=\left\{ Ap\ :\ A\in \Pi_{k-m}\right\},\quad
{\mathcal W_2}=\left\{ Bq\ :\ B\in \Pi_{k-n}\right\}.$$
Since $p$ and $q$ have no common component we conclude that
$${\mathcal W_1}\cap{\mathcal W_2}=\left\{ Cpq\ :\ C\in \Pi_{k-m-n}\right\}.$$
Now we readily obtain that
\begin{multline}\label{er}\dim{\mathcal W}=\dim({\mathcal W_1}+{\mathcal W_2})=\dim{\mathcal W_1}+\dim{\mathcal W_2}-\dim({\mathcal W_1}\cap{\mathcal W_2})\\={k-m+2\choose 2}+{k-n+2\choose 2}-{k-m-n+2\choose 2}\\={k+2\choose 2}-mn.
\end{multline}
The last equality here holds since $k\ge m+n$ (actually it holds for $k\ge m+n-2$).

{\bf Step 2.} $ n+m\le k\le k_0.$

Let us apply decreasing induction with respect to $k.$ The first step $k=k_0$ was checked in Step 1.
Assume Theorem is true for all $f$ with $\deg f=k$ and let us prove that it is true also for all $f$ with $\deg f=k-1.$

Suppose that $f_0$ is an arbitrary polynomial with $\deg f_0=k-1.$ Choose a line $\ell_0$
such that

(i) $\ell_0\cap p\cap q=\emptyset,$ and

(ii) $\ell_0$ intersects $q$ at $n$ points, counted also multiplicities, i.e., it does not intersect
$q$ at infinity.

We have that $\deg f_0\ell_0=k.$ Also, in view of the formula \eqref{main} and \eqref{inv}, i.e., the $D$-invariance of $\Mset_\lambda(p,q)$, we have that $f_0\ell_0$ vanishes at $\Mset_\lambda(p,q)$ for each $\lambda\in p\cap q.$ Hence, in view of the induction hypothesis, we get
\begin{equation}\label{fl} f_0\ell_0=Ap+Bq,
\end{equation}
where $A\in  \Pi_{k-m},\ B\in \Pi_{k-n}.$

We have that $\ell_0$ intersects $q$ at $n$ points, counted also multiplicities. In view of \eqref{fl} these (multiple) points are also zeros of $A$ since $p$ differs from zero there.

For every polynomial  $C_0\in\Pi_{k-m-n}$ we have also that
\begin{equation}\label{fl2}f_0\ell_0=(A-C_0q)p+(B+C_0p)q.\end{equation}
Consider arbitrary $k-m-n+1$ points $\lambda_1,\ldots,\lambda_{k-m-n},$ in $\ell_0\setminus q.$ Choose $C_0\in\Pi_{k-m-n}$ such that $A-C_0q$ is zero at these points.
For this, according to Proposition \ref{n+1points}, we just solve an independent interpolation problem
$$ C_0(\lambda_i)=\frac {A(\lambda_i)}{q(\lambda_i)},\quad i=0,\ldots,k-m-n.$$

Note that the common $n$ (multiple) zeros of $\ell_0$ and $q$ also are zeroes of $A-C_0q.$ Thus, altogether we have that  $A-C_0q$ is zero at $k-m-n+1+n=k-m+1$ points in $\ell_0.$
Thus, in view of Proposition \ref{n+1points}, $\ell_0$ divides $A-C_0q\in \Pi_{k-m}.$ From \eqref{fl2} we readily conclude that $\ell_0$ divides $B+C_0p.$

Finally by dividing the relation \eqref{fl2} by $\ell_0$ we get that
\begin{equation}\label{fl3} f_0\ell_0=A'p+B'q,
\end{equation}
where $A'\in  \Pi_{k-m-1},\ B\in \Pi_{k-n-1}.$

{\bf Step 3.} $k\le n+m-1.$

Let us again apply decreasing induction with respect to $k.$ The first step $k=m+n-1$ was checked in Step 2.
Assume Theorem is true for all $f$ with $\deg f=k$ and let us prove that it is true also for all $f$ with $\deg f=k-1.$

Suppose that $f_0$ is an arbitrary polynomial with $\deg f_0=k-1.$ Choose a line $\ell_0$ in the same way as in Step 2.
Then we get the relation \eqref{fl} where the polynomial $A\in  \Pi_{k-m}$ has $n$ zeros in $\ell_0,$ counting also the multiplicities.
In this case we have that $k-m\le n-1.$
Thus, in view of Proposition\ref{n+1points}, $\ell_0$ divides $A.$ From \eqref{fl2} we readily conclude that $\ell_0$ divides also $B.$
Finally by dividing the relation \eqref{fl} by $\ell_0$ we complete the proof as in Step 2.
\hfill$\Box$

At the end let us bring the formulation of Theorem \ref{Noether} in the homogeneous case.

\begin{theorem} \label{Be} Suppose that $p\in \Pi_m^0$ and $q\in\Pi_n^0$ have no common component. Suppose also that $f\in \Pi_k^0$ vanishes at $\Mset_\lambda^0(p,q)$ for each $\lambda\in p\cap q.$ Then we have that
\begin{equation*}f=Ap+Bq,
\end{equation*}
where $A\in \Pi_{k-m}^0,\ B\in \Pi_{k-n}^0.$
\end{theorem}

It is known that the set $\Zset_0:=p\cap q,$ where $p$ and $q$ are polynomials,  of degree  $m$ and $n,$ respectively, is $(m+n-2)$-independent, provided that $|\Zset_0|=mn.$
Below we prove this result without the last restriction (cf. \cite{H04}, Corollary 1).
\begin{corollary} \label{AA} Suppose that polynomials $p,q\in \Pi,\ \deg p=m,\ \deg q=n,$ have no common component. Then the set of linear operators $\Lset(p,q)$ and consequently the set $\Zset_0$ are $(m+n-2)$-independent.
\end{corollary}
\textbf{Proof.}
Let us assume first that $p$ and $q$ have no intersection point at infinity.  Then we have that $|\Lset(p,q)|=mn.$
By using the evaluation \eqref{er} in the case $k=m+n-2$ we obtain
\begin{multline}\dim{\mathcal W}=\dim({\mathcal W_1}+{\mathcal W_2})=\dim{\mathcal W_1}+\dim{\mathcal W_2}-\dim({\mathcal W_1}\cap{\mathcal W_2})\\={n\choose 2}+{m\choose 2}-0={m+n\choose 2}-mn.
\end{multline}
Thus we have that $\dim\Pi_{m+n-2}-\dim{\mathcal W}=mn.$ This means that
the set of linear operators $\Lset(p,q)$ and consequently $\Zset_0$ is $(m+n-2)$-independent.

Now assume only that $p$ and $q$ have no common component.
Let us use the concept of the associate polynomial (see section 10.2, \cite{W}). 

Let $p(x,y)=\sum_{i+j\le m}a_{ij}x^iy^j$ and $\deg p=m.$ Then the following trivariate homogeneous polynomial is called associated with $p:$
$$\bar p(x,y,z)= \sum_{i+j+k=m}a_{ij}x^iy^jz^{k}.$$
Evidently we have that
$$
p=p_1p_2 \Leftrightarrow \bar p=\bar p_1\bar p_2.
$$
It is easily seen from here that polynomials $p$ and $q$ have no common component if and only if $\bar p$ and $\bar q$ have no common component.
By applying Theorem \ref{Be} to the polynomials $\bar p$ and $\bar q$ we get that
the set of linear operators $\Lset^0(p,q)$ is $(m+n-2)$-independent. Therefore its subset
corresponding to the finite intersection points, i.e., to $\Zset_0,$ is $(m+n-2)$-independent, which implies the desired result.
\hfill$\Box$

\section {The Cayley-Bacharach theorem}

The evaluation \eqref{er} in the case $k=m+n-3$ gives
\begin{multline}\dim{\mathcal W}=\dim({\mathcal W_1}+{\mathcal W_2})=\dim{\mathcal W_1}+\dim{\mathcal W_2}-\dim({\mathcal W_1}\cap{\mathcal W_2})\\={n-1\choose 2}+{m-1\choose 2}-0={m+n-1\choose 2}-(mn-1).
\end{multline}
Thus we have that $\dim\Pi_{m+n-2}-\dim{\mathcal W}=mn-1,$ i.e., out of $mn$ linear operators in $\Lset(p,q)$ only $mn-1$ are linearly independent.

According to the Cayley-Bacharach classic theorem (see, e.g., \cite{EGH}, \cite{HJZ}), i.e., in the case $|{\mathcal Z}_0| = mn,$ where ${\mathcal Z}_0:=p\cap q,$ we have that any    subset of $\Zset_0$ of cardinality $mn-1$ is $(m+n-3)$-independent.
This means that no point from $\Zset_0$ has a fundamental polynomial of degree $(m+n-3),$ i.e., for any point $\lambda_0\in \Zset_0$ the following implication holds:
$$ p\in\Pi_{m+n-3},\ p(\lambda)=0 \ \forall \lambda\in \Zset_0\setminus \{\lambda_0\}
\Rightarrow p(\lambda)=0 \ \forall \lambda\in \Zset_0.$$

In this section we are going to study the situation in the general multiple intersection case.
Suppose $p\in \Pi_m,$
 \begin{equation*}
p(x,y)=\sum_{i+j\le m}a_{ij}x^iy^j.
\end{equation*}
Denote the $k$th homogeneous part of $p$ by $p^{\{k\}},$ i.e.,
\begin{equation*}
p^{\{k\}}(x,y)=\sum_{i+j=k}a_{ij}x^iy^j.
\end{equation*}
We accept a very common restriction from the theory of intersection. Namely, we assume that the two polynomials $p$ and $q$ have no common tangent line at an intersection point $\lambda\in\Zset_0.$ This means that the lowest homogeneous parts of the polynomials have no common factor at this point.

\begin{theorem} Suppose that polynomials $p,q\in \Pi,\ \deg p=m,\ \deg q=n,$ have no intersection point at infinity and $\lambda\in \Zset_0.$ Suppose also that $p$ and $q$ have no common tangent line at $\lambda.$ Then we have that the set of linear operators $\Lset^\lambda(p,q)$ contains only one operator of the highest degree: $\bar L.$ Suppose also that $f\in  \Pi_{m+n-3}$ vanishes at $\Lset(p,q)\setminus \{\bar L\}.$ Then we have that $f$ vanishes at all $\Lset(p,q).$
\end{theorem}
\textbf{Proof.} Assume, without loss of generality, that $\lambda=\theta:=(0,0).$
Suppose that $p$ and $q$ are bivariate polynomials having  $n_0$ and $m_0$-fold zero at the origin, respectively, $n_0,m_0\ge 1$ :
\begin{equation*}
p(x,y)=\sum_{m_0\le i+j\le m}a_{ij}x^iy^j, \quad q(x,y)=\sum_{n_0\le i+j\le n}b_{ij}x^iy^j.
\end{equation*}
Suppose also that  $p$ and $q$  have no common tangent line at the origin, i.e., $p^{\{m_0\}}$ and $q^{\{n_0\}}$ have no common factor.

Let $\bar\Lset:=\{\bar L_1,\ldots,\bar L_s\}$ be a maximal independent set of linear  operators with the highest degree  in the  space $\Mset_\theta(p,q).$

Assume that $f\in  \Pi_{m+n-3}$ vanishes at $\Lset(p,q)\setminus \bar \Lset.$
We are going to prove that $f$ vanishes at $\Lset(p,q).$

This shall complete the proof of Theorem. Indeed, as was verified above, there are   $mn-1$ linearly independent operators in the set of $mn$ linear operators $\Lset(p,q),$ which clearly implies here that $s=1.$

Let $\ell$ be any line passing through $\theta.$
By using the formula \eqref{main} with $g=\ell,\ f=f$ and $R\in \Lset(p,q),$ we obtain that the polynomial $\ell f$ vanishes  at $\Lset(p,q).$
Therefore, since $\deg \ell f =m+n-2,$ we get from Theorem \ref{Noether} that
\begin{equation}\label{aaalf}\ell f=A(\ell)p+B(\ell)q,
\end{equation}
where $A(\ell)\in \Pi_{n-2},\ B(\ell)\in \Pi_{m-2}.$
Assume, without loss of generality, that $m_0\le n_0.$
Assume also that $m_0\ge 2.$ If $m_0=1$ we go to the \emph{final part of the proof}.
Now we are going to prove that
\begin{equation}\label{Al} A(\ell)^{\{k\}}=\ell A'_{k-1}\quad k=0,\ldots,n_0-2,
\end{equation}
where $A'_{k-1}, \in \Pi_{k-1}^0,$ do not depend on $\ell,$ and
\begin{equation}\label{Bl} B(\ell)^{\{k\}}=\ell B'_{k-1}\quad k=0,\ldots,m_0-2,
\end{equation}
where $B'_{k-1}, \in \Pi_{k-1}^0,$ do not depend on $\ell.$

First let us prove \eqref{Al} for $k\le n_0-m_0-1.$ Let us apply induction on $k.$
Consider the case $k=0.$
Then we get from the relation \eqref{aaalf} that
$A(\ell)^{\{0\}}p^{\{m_0\}}=\ell f^{\{m_0-1\}}.$
Thus we have
$xf^{\{m_0-1\}}=c_1p^{\{m_0\}}$ and $yf^{\{m_0-1\}}=c_2p^{\{m_0\}},$
where $c_1$ and $c_2$ are constants. Therefore we have that $(c_2x-c_1y)f^{\{m_0-1\}}=0,$ i.e., $f^{\{m_0-1\}}=0.$ Thus $A(\ell)^{\{0\}}=0=\ell \cdot0.$
Assume that \eqref{Al} is true for all $k$ not exceeding $s$ and let us prove it for $k=s+1.$
We readily get from the relation \eqref{aaalf} that
\begin{equation}\label{gum}A(\ell)^{\{s+1\}}p^{\{m_0\}}+A(\ell)^{\{s\}}p^{\{m_0+1\}}+\cdots +A(\ell)^{\{0\}}p^{\{m_0+s+1\}}=\ell f^{\{m_0+s+1\}}.\end{equation}
We have that all terms above except possibly the first have factor $\ell.$ Hence we get
that $A(\ell)^{\{s+1\}}=\ell A'_s.$ In fact we have this relation for all $\ell$ except $m_0$ tangent lines of $p$ at $\theta.$ Then by a continuity argument we get the relation for all $\ell.$

Next, by dividing \eqref{gum} by $\ell$ we see that $A'_s$ does not depend on $\ell.$

Now assume that $n_0-m_0\le k\le n_0-2.$ Here we are going to prove \eqref{Al} for
$k$ and \eqref{Bl} for $k-n_0+m_0.$ Let us again apply induction on $k.$
Consider the case $k=n_0-m_0.$
We get from the relation \eqref{aaalf} that
\begin{multline}\label{gum2}A(\ell)^{\{n_0-m_0\}}p^{\{m_0\}}+A(\ell)^{\{n_0-m_0-1\}}p^{\{m_0+1\}}+\cdots +A(\ell)^{\{0\}}p^{\{n_0\}}+B(\ell)^{\{0\}}q^{\{n_0\}}\\=\ell f^{\{n_0-1\}}.\end{multline}

Now let us use $\ell=\ell_1$ which is a tangent line of $q$ at $\theta,$ i.e.,  $q^{\{n_0\}}=\ell_1 \tilde q,$ where $\tilde q\in \Pi_{n_0-1}.$

We have that all terms in \eqref{gum2} except possibly the first have factor $\ell_1.$ Hence we get
that $A:=A(\ell_1)^{\{n_0-m_0\}}=\ell_1 A'_{n_0-m_0-1}.$ 

Meanwhile, let us verify also that if $\ell_1=y-k_1x$ is a factor of multiplicity $\mu$ of  $q^{\{n_0\}}$ then it is a factor of  multiplicity at least $\mu$ in $A.$
Assume that $$A=C_1\prod_{i}(y-a_ix), \ q^{\{n_0\}}=C_2\prod_{i}(y-b_ix).$$
Assume also $\ell$ is given by an equation $y-kx=0.$
By setting in \eqref{gum2} $y=kx,$ and by using the induction hypothesis, we obtain
\begin{equation}\label{pq}C_1p^{\{m_0\}}(x,kx)\prod_{i}(k-a_i)x=C_2B(\ell)^{\{0\}}(x,kx)\prod_{i}(k-b_i)x.
\end{equation}
Consider both sides of \eqref{pq} as polynomials on $k.$
Now $k_1$ is a root of the right hand side of multiplicity at least $\mu.$ On the other hand $k=k_1$
is not a root of $p^{\{m_0\}}(x,kx)$ since $p$ and $q$ have no common factor.
Thus we get that $k=k_1$ is a root of multiplicity at least $\mu$ in $q^{\{n_0\}}(x,kx),$
i.e., $y-k_1x$ is a factor of multiplicity at least $\mu$ in $q^{\{n_0\}}(x,y).$

Next, we have that \begin{multline}A(\ell)^{\{n_0-m_0\}}=A(\ell_1)^{\{n_0-m_0\}}+A(\ell-\ell_1)^{\{n_0-m_0\}}\\=\ell A'_{n_0-m_0-1}+(\ell-\ell_1) A'_{n_0-m_0-1}+A(\ell-\ell_1)^{\{n_0-m_0\}}\\=\ell A'_{n_0-m_0-1}-(k-k_1)x A'_{n_0-m_0-1}-(k-k_1)A(x)^{\{n_0-m_0\}}\\
=\ell A'_{n_0-m_0-1}-(k-k_1)\left[x A'_{n_0-m_0-1}-A(x)^{\{n_0-m_0\}}\right].\end{multline}
We have that  $A(\ell)^{\{n_0-m_0\}}$ contains all factors of $q_{n_0}.$ Thus the polynomial of degree $n_0-m_0$ in the square brackets contains all factors of $q_{n_0}$ except possibly $\ell_1,$ in all $n_0-1$ factors.
Hence this polynomial is identically zero and
$A(\ell)^{\{n_0-m_0\}}=\ell A'_{n_0-m_0-1}.$
As above we readily conclude that $A'_{n_0-m_0-1}$ does not depent on $\ell.$
Similarly by using tangent lines of $p$ we get that $B(\ell)^{\{0\}}=0=\ell\cdot 0.$

Now assume that \eqref{Al} is true for $k$ not exceeding $s$ and \eqref{Bl} is true for $k$ not exceeding $s+m_0-n_0.$ Let us prove \eqref{Al} for $k=s+1$ and \eqref{Bl}
for $k=s+m_0-n_0+1.$

We get from the relation \eqref{aaalf} that
\begin{multline}\label{gum}A(\ell)^{\{s+1\}}p^{\{m_0\}}+A(\ell)^{\{s\}}p^{\{m_0+1\}}+\cdots +A(\ell)^{\{0\}}p^{\{m_0+s+1\}}\\+B(\ell)^{\{s+m_0-n_0+1\}}q^{\{n_0\}}+B(\ell)^{\{s+m_0-n_0\}}q^{\{n_0+1\}}+\cdots +B(\ell)^{\{0\}}q^{\{m_0+s+1\}}\\=\ell f^{\{m_0+s+1\}}.\end{multline}
Here, in the same way as above, by using tangent lines of $p$ and $q$ at $\theta$, we
complete the proof of this part.

Now let us go to the \emph{final part of the proof}.
Let us choose a line $\ell_0$ whose intersection multiplicity with $p$ at $\theta$
equals to $m_0.$ We also require that $\ell_0$ intersects $\Zset$ only at $\theta.$
We have that outside of $\theta$ the line $\ell_0$ intersects $p$ at $m-m_0$ points, counting also the multiplicities.
We deduce from the relation \eqref{aaalf}, with $\ell=\ell_0,$ that these $m-m_0$ points
are roots for $B(\ell_0),$ since $q$ does not vanish there.
Then, in view of the relation \eqref{Bl}, we have that
$$B(\ell_0)=\sum_{i=0}^{m-2}B^{\{i\}}(\ell_0)=\sum_{i=m_0-1}^{m-2}B^{\{i\}}(\ell_0).$$
Thus, by assuming that $\ell_0=y-k_0x,$ we see that the trace of the polynomial $B(\ell_0)$ on the line $\ell_0$ has the form $$B(\ell_0)(x,k_0x)=\sum_{i=m_0-1}^{m-2}b_ix^i=x^{m_0-1}\sum_{i=0}^{m-m_0-1}b_{i+m_0-1}x^i.$$ On the other hand this polynomial
vanishes at $m-m_0$ nonzero points, counting also the multiplicities. Hence, in view of  Proposition \ref{n+1points}, we conclude that
$B(\ell_0)$ has a factor $\ell_0.$ Now we readily get from the relation \eqref{aaalf}, with $\ell=\ell_0,$ that $A(\ell_0)$ also has a factor $\ell_0.$ Then by dividing the relation \eqref{aaalf} by $\ell_0$ we get that
\begin{equation*}f=Ap+Bq,
\end{equation*}
where $A\in \Pi_{n-3},\ B\in \Pi_{m-3}.$
Finally from this relation we readily conclude that $f$ vanishes at $\Lset(p,q).$
\hfill$\Box$

At the end let us consider a simple example. Let $p(x,y)=x^m$ and $q(x,y)=y^n.$
Then we have that $$\Lset(p,q)=\Lset_\theta(p,q)=\left\{x^iy^j : i\le m-1,\ j\le n-1\right\}.$$
It is easily seen that in this set there is only one operator of the highest degree: $$\bar L=
\left(\frac{\partial }{\partial x}\right)^{m-1}\left(\frac{\partial }{\partial y}\right)^{n-1}.$$ Also for this operator we have that the set of the operators $\Lset(p,q)\setminus \{\bar L\}$ is $(m+n-3)$-independent. Moreover, only the operator $\bar L\in \Lset(p,q)$ has this property.

\end{document}